\documentclass[a4paper,12pt]{article}
\usepackage[latin1]{inputenc}
\usepackage{amsmath,amssymb,amsthm,cite,verbatim}
\usepackage[all]{xy}
\usepackage[bookmarks=true]{hyperref}
\usepackage{color}
\usepackage{graphicx}

\title{The automorphism group of the non-split Cartan modular curve of level 11}
\author{Valerio Dose, \ Julio Fern\'andez, \ Josep Gonz\'alez  \ and \ Ren\'e Schoof \footnote{The second and the third authors are partially supported by DGI grant  MTM2012-34611.
}}

\newtheorem{prop}{Proposition}
\newtheorem{cor}{Corollary}
\newtheorem{rem}{Remark}
\newtheorem*{xthm}{Theorem}

\theoremstyle{definition}

\theoremstyle{remark}

\numberwithin{equation}{section}

\newcommand{\QQ}{\mathbf{Q}}

\newcommand{\QQbar}{\overline{\mathbf{Q}}}

\newcommand{\ZZ}{\mathbf{Z}}

\newcommand{\FF}{\mathbf{F}}

\newcommand{\PP}{\mathbf{P}}

\newcommand{\CC}{\mathbf{C}}

\newcommand{\RR}{\mathbf{R}}

\newcommand{\Gal}{\mathrm{Gal}}

\begin{document}
\maketitle

\begin{abstract}
\noindent
We derive equations for the modular curve $X_{ns}(11)$ associated to  a non-split Cartan subgroup of ${\rm GL}_2(\FF_{11})$. This allows us to compute the automorphism group of the curve and show that it is isomorphic to Klein's four group.
\end{abstract}

\section*{\large Introduction}

Let $p$ be a prime. The modular curve $X_{ns}(p)$ associated to a non-split Cartan subgroup of ${\rm GL}_2(\FF_p)$
is an algebraic curve that is defined over~$\QQ$. It admits a so-called {\it modular involution}~\,$w$,  also
defined over~$\QQ$. One may conjecture that, for large $p$, the modular involution is the only non-identity
automorphism of $X_{ns}(p)$, even over $\CC$.  However, for very small primes $p$ this is not the case.
Indeed, for $p=2,3$ and $5$ the genus of $X_{ns}(p)$ is~$0$,
while for $p=7$ the genus is $1$. See~\cite{BAR}, Table~A.1.
For these primes the curve $X_{ns}(p)$ admits therefore  infinitely many automorphisms. The present paper is
devoted to $p=11$ and the genus $4$ curve $X_{ns}(11)$. We prove the following.

\begin{xthm}\label{theorem}
The automorphism group over $\CC$ of the modular curve $X_{ns}(11)$
is isomorphic to Klein's four group. It is generated by the modular involution $w$ and
the  involution $\varrho$ described in Corollary~\rm{\ref{exotic}}.
\end{xthm}

Our proof for this result is presented in section~\ref{aut}.
It relies on an explicit description of the regular differentials
and the Jacobian of~$X_{ns}(11)$. These  are discussed in section~\ref{dif}.
We make use of equations for the curve $X_{ns}(11)$, which are obtained in section~\ref{equ}.

\section{\large Equations}\label{equ}

In this section we derive equations for the modular curve $X_{ns}(11)$. We do this by exploiting the
modular curve $X_{ns}^+(11)$ associated to the normalizer of a non-split Cartan subgroup of level~$11$.

We recall some definitions~\cite{BAR}.
For any prime $p$, the ring of $2\times 2$ matrices over~$\FF_p$ contains subfields that are isomorphic to $\FF_{p^2}$.
A non-split Cartan subgroup $U$ of ${\rm GL}_2(\FF_p)$ is by definition the unit group of such a subfield.
The modular curve $X_{ns}(p)$ classifies $U$--isomorphism classes of pairs $(E,\phi)$, where $E$ is an elliptic
curve and $\phi$ is an isomorphism from the group of $p$-torsion points $E[p]$
to $\FF_p\!\times\!\FF_p$. Two such pairs $(E,\phi)$ and $(E',\phi')$ are $U$--isomorphic if there is
an isomorphism $f:E\longrightarrow E'$ for which the matrix $\phi'f\phi^{-1}$ is in $U$.

The group $U$ has index $2$ in its normalizer $U^+\!\subset {\rm GL}_2(\FF_p)$.
The modular involution $w$ of $X_{ns}(p)$ maps $(E,\phi)$ to $(E,\alpha\phi)$,
where $\alpha$ is any matrix in $U^+\!\smallsetminus\!U$.
In a way that is analogous to the moduli description for $X_{ns}(p)$, the modular curve $X_{ns}^+(p)$
classifies $U^+$--isomorphism classes of pairs~$(E,\phi)$. There are natural morphisms
$$
X_{ns}(p)\mathop{\longrightarrow}\limits^{\pi} X_{ns}^+(p)\mathop{\longrightarrow}\limits^{j} X(1).
$$
Here $X(1)$ indicates the $j$--line. It parametrizes elliptic curves up to isomorphism.
The morphism $j$ maps $(E,\phi)$ to the $j$--invariant of $E$. It has degree ${1\over 2}\,p(p-1)$,
while the morphism $\pi$ has degree~$2$.

Both curves $X_{ns}(p)$ and $X_{ns}^+(p)$ are defined over $\QQ$. A point $(E,\phi)$ of
$X_{ns}(p)$ or $X_{ns}^+(p)$ is defined over an extension $\QQ\subset K$
if and only if  $E$ is defined over $K$ and, for all $\sigma\in{\rm Gal}(\overline{K}/K)$,
the matrix $\phi\sigma\phi^{-1}$ is in $U$ or $U^+$ respectively.
This implies that, for $p>2$, the curve $X_{ns}(p)$ does not contain any points defined over $\RR$.
On the other hand, the curve $X_{ns}^+(p)$ has real and usually also rational points.
Indeed, for every imaginary quadratic order~$R$ with class number~$1$ there is a unique elliptic curve $E$ over $\CC$
with complex multiplication by $R$. The $j$--invariant of~$E$ is in~$\QQ$. Moreover, when $p$ is  prime in the ring $R$,
there is a unique rational point $(E,\phi)$ on~$X_{ns}^+(p)$.
These points are called {\it CM points} or {\it Heegner points}. See~\cite{APP},~Section A.5.

\begin{rem}\label{quadratic}{\rm
Suppose that $(E,\phi)$ is a rational point of $X_{ns}^+(p)$. Then $E$ is defined over~$\QQ$ and
the image of ${\rm Gal}(\overline{\QQ}/\QQ)$ in ${\rm Aut}(E[p])$ is isomorphic through $\phi$ to a subgroup $G$ of
${\rm GL}_2(\FF_p)$ which is contained in the normalizer of a non-split Cartan subgroup $U$.
The points of $X_{ns}(p)$ lying above $(E,\phi)$ are defined over the fixed field of \,$U\cap G$,
which is an imaginary quadratic extension of~$\QQ$. In the case of Heegner points, CM theory implies
that this extension is isomorphic to the quotient field of the endomorphism ring of~$E$.}
\end{rem}

Now we turn to the case $p=11$.
In \cite{LIG},~Proposition 4.3.8.1, Ligozat derived a Weierstrass equation for the genus $1$ curve $X_{ns}^+(11)$. It is given by
$$
Y^2+Y=X^3-X^2-7X+10.
$$
By choosing the point at infinity as origin, we can view $X_{ns}^+(11)$ as an elliptic curve
and equip it with the usual group law.  The rational points of $X_{ns}^+(11)$ are then an infinite cyclic group
generated by the point $P=(4,-6)$.  See~\cite{CRE}.
The  translations by the rational points form an  infinite  group of  automorphisms
of $X_{ns}^+(11)$. They are all  defined over~$\QQ$.
It follows that there are infinitely many isomorphisms over $\QQ$
between $X_{ns}^+(11)$ and the curve given by Ligozat.  For a particular choice of such an isomorphism,
Halberstadt derived in \cite{HALB}, Section~2.2, an explicit formula for the degree $55$ morphism
\,$j\!:X_{ns}^+(11)\longrightarrow X(1)$. In view of the symmetry phenomenon described at the end of this section,  it is convenient to compose his isomorphism with the
translation-by-$P$ morphism. Explicitly, our function $j(X,Y)$ is  the value of Halberstadt's \,$j$--function in the point
$$\left({{4X^2+X-2+11Y}\over{(X-4)^2}}\,,\,{{(2X^2+17X-34+11Y)(1-3X)}\over{(X-4)^3}}\right)\!,$$
that is,
$$
\begin{array}{l}
\!\!\!\!j(X,Y)=(X+2)(4-X)^5\,\big(11(X^2+3X-6)(Y-5)(X^3+4X^2+X+22+(1-3X)Y)\big)^3\\[8pt]
\,\times\,\dfrac{\big((3X^2-3X-14-(3+2X)Y)(12X^3+28X^2-41X-62+(3X^2+20X+37)Y)\big)^3}
{\big(-7X^2-15X+62+(X+18)Y\big)^2\,\big(4X^3+2X^2-21X-6+(X^2+3X+5)Y\big)^{11}}\,.
\end{array}
$$

\begin{prop}\label{equations}  The modular curve $X_{ns}(11)$ is given by the equations
$$
\begin{array}{rcl}
Y^2+Y & \!=\! & X^3-X^2-7X+10,\\[5pt]
T^2 & \!=\! & -(4X^3+7X^2-6X+19).
\end{array}$$
\end{prop}

\noindent{\bf Proof.}\,  We first compute the ramification locus of the morphism
$\pi\,\colon X_{ns}(11)\!\longrightarrow\!X^+_{ns}(11)$. Since~$\pi$ is defined over $\QQ$,
this locus is stable by the action of~\,$\Gal(\QQbar/\QQ)$.
By Proposition~7.10 in \cite{BAR}, the function \,$j(X,Y)\!-\!1728$\, has exactly seven simple zeroes on $X^+_{ns}(11)$,
and six of them are the rami\-fi\-ca\-tion points of $\pi$. All the other zeroes are double.
Let us consider the quotient map \,$X_{ns}^+(11)\longrightarrow \PP^1$ induced by the elliptic involution.
It corresponds to the quadratic function field \mbox{extension \,$\QQ(X)\subset\QQ(X,Y)$} with non-trivial
automorphism given by \,$Y\mapsto -1-Y$\!. One easily checks that the trace and norm of the function
$j(X,Y)\!-\!1728$ admit the polynomial \,$4X^3+7X^2-6X+19$\, as an irreducible factor
of multiplicity $1$ and $2$ respectively. The function~$F$ on $X^+_{ns}(11)$ defined by this cubic
polynomial has exactly six simple zeroes. It follows that the zeroes of $F$ are simple zeroes of
\,$j(X,Y)\!-\!1728$. Therefore they are the ramification points of~$\pi$.

The function field $\QQ(X_{ns}(11))$ is obtained by adjoining a function $G$ to $\QQ(X_{ns}^+(11))$ whose square
is in $\QQ(X_{ns}^+(11))$. The coefficients of the divisor on $X_{ns}^+(11)$ of \,$G^2$ are odd at the ramified points
and even at the others. Since the same holds for the above function $F$\!,\, the divisor of~$F\,G^2$ is of the form
\,$2D$\, for some divisor \,$D$\, of $X_{ns}^+(11)$ defined over~$\QQ$.
The group ${\rm Pic}^0(X_{ns}^+(11))$ is naturally isomorphic to the group of rational points of $X_{ns}^+(11)$.
Since the  latter is isomorphic to~$\ZZ$, there are no elements of order~$2$ in~${\rm Pic}^0(X_{ns}^+(11))$.
It follows that $D$ is principal.
This means that there is a function \,$T$\, in $\QQ(X_{ns}(11))$ and a non-zero \,$\lambda\!\in\!\QQ$\, for which
$\lambda\,T^2=F$\!.\, The function field of~$X_{ns}(11)$ is then equal to~$\QQ(X,Y,T)$.

It remains to determine  $\lambda$, which is unique up to squares. Consider the point \,$Q=(5/4,7/8)$ of $X_{ns}^+(11)$.
Since $j(Q)=1728$, the elliptic curve parametrized by the point $Q$ admits complex multiplication by the ring $\ZZ[i]$
of Gaussian integers. By Remark~\ref{quadratic}, the two points of $X_{ns}(11)$ lying above $Q$ are defined over $\QQ(i)$. Since \,$F(Q)\!=\!121/4$\, is a square, we may take~$\lambda=-1$. This proves the proposition.

\vspace{1truemm}

\begin{cor}\label{exotic}
In addition to the modular involution $w$, the curve $X_{ns}(11)$ admits
an \mbox{``exotic"} involution $\varrho$. The modular involution  switches  $(X,Y,T)$ and\, $(X,Y,-T)$, while
$\varrho$ switches $(X,Y,T)$ and\, $(X,-1-Y,T)$. Together, $w$ and $\varrho$ generate a  subgroup
of~\,${\rm Aut}(X_{ns}(11))$  iso\-morphic to Klein's four group.
\end{cor}

Although it is not relevant for the proofs in this paper, let us explain how the ``exotic"
automorphisms of $X_{ns}(11)$ were first detected.  The rational points of $X_{ns}^+(11)$ form an infinite
cyclic group generated by the point $P=(4,-6)$. For each $n\in\ZZ$, the elliptic curve over $\QQ$ parametrized by
the point $[n]P$ in $X_{ns}^+(11)(\QQ)$ has the following property: the image $G$ of the Galois representation
attached to its $p$-torsion points is contained in the normalizer of a non-split Cartan subgroup $U$.
By Remark~\ref{quadratic}, the fixed field of \,$U\cap G$ is an imaginary quadratic field.
In his {\it tesi di laurea}, one of the authors
---Valerio Dose--- used the methods of~\cite{JPS} to compute this quadratic field $K$ for several values of $n$.
The first few values are given in the table below. There is a striking symmetry: the quadratic fields
attached to the points $[n]P$ and $[-n]P$ are always the same. There does not seem to be a ``modular reason" for this,
as it may happen that the elliptic curve associated to $[n]P$ has complex multiplication by some quadratic order of
discriminant $\Delta<0$ but the elliptic curve associated to $[-n]P$ has not. In the first case $K$ is the CM field,
but in the second case it is not. The phenomenon, which surprised us at first, is explained by the existence of the
``exotic" involution~$\varrho$.
\begin{center}
\begin{tabular}{|r|cc|cl|cc|}
\hline
&&&&&&\\[-10pt]
points && $j$ && \quad CM && $K$\\[5pt]
\hline
&&&&&&\\[-10pt]
$[6]P$ \ &&$2^33^95^311^317^629^353^3191^3/769^{11}$&&\quad \ --&&$\QQ(\sqrt{-3\cdot 14327})$\\[5pt]
$[5]P$ \ &&$-2^183^35^323^329^3$&&$\Delta=-163$&&$\QQ(\sqrt{-163})$\\[5pt]
$[4]P$ \ &&$0$&&$\Delta=-3$&&$\QQ(\sqrt{-3})$\\[5pt]
$[3]P$ \ &&$2^63^3$&&$\Delta=-4$&&$\QQ(\sqrt{-1})$\\[5pt]
$[2]P$ \ &&$-2^{15}3^35^311^3$&&$\Delta=-67$&&$\QQ(\sqrt{-67})$\\[5pt]
$P$ \ &&$2^43^35^3$&&$\Delta=-12$&&$\QQ(\sqrt{-3})$\\[5pt]
$\infty$ \ &&$2^33^311^3$&&$\Delta=-16$&&$\QQ(\sqrt{-1})$\\[5pt]
$[-1]P$ \ &&$-2^{15}3^15^3$&&$\Delta=-27$&&$\QQ(\sqrt{-3})$\\[5pt]
$[-2]P$ \ &&$2^83^35^611^353^3/23^{11}$&&\quad \ --&&$\QQ(\sqrt{-67})$\\[5pt]
$[-3]P$ \ &&$-2^93^35^313^171^3181^3/43^{11}$&&\quad \ --&&$\QQ(\sqrt{-3})$\\[5pt]
$[-4]P$ \ &&$2^{18}3^35^37^111^323^329^3103^3/67^{11}$&&\quad \ --&&$\QQ(\sqrt{-3})$\\[5pt]
$[-5]P$ \ &&$-2^43^35^117^629^3367^32381^3/397^{11}$&&\quad \ --&&$\QQ(\sqrt{-163})$\\[5pt]
$[-6]P$ \ &&$-2^33^111^317^619^123^341^353^3167^32777^323431^3/80233^{11}$&&\quad \ --&&$\QQ(\sqrt{-3\cdot 14327})$\\[5pt]
\hline
\end{tabular}
\end{center}

\vspace{1truemm}

\section{\large Differentials}\label{dif}

In this section we analyze the space of regular differentials \,$\Omega^1_{X_{ns}(11)}$ of the curve $X_{ns}(11)$.

By \cite{CHEN}, Section~8, the Jacobian \,$J_{ns}(11)$ of $X_{ns}(11)$ is isogenous over $\QQ$ to the new part
of the Jacobian of $X_0(121)$. See \cite{EDIX} for an easy proof of this result. By Cremona's Tables~\cite{CRE},
there are exactly four $\QQ$--isogeny classes of elliptic curves of conductor $121$, which are represented by
$$
\begin{array}{ccl}
A&:&\quad y^2+xy+y=x^3+x^2-30x-76,\\[6pt]
B&:&\quad y^2+y=x^3-x^2-7x+10,\\[6pt]
C&:&\quad y^2+xy=x^3+x^2-2x-7,\\[6pt]
D&:&\quad y^2+y=x^3-x^2-40x-221.
\end{array}
$$
It follows that $J_{ns}(11)$ is isogenous over $\QQ$ to the product of these four elliptic curves.
The following proposition describes a low degree morphism from the curve~$X_{ns}(11)$ to each of its elliptic
quotients, and provides a basis for \,$\Omega^1_{X_{ns}(11)}$ from the respective pull-backs.
We make use of the equations for $X_{ns}(11)$ given in Proposition \ref{equations}.
It is also convenient to introduce the function \,$Z=(2Y+1)T$\, in \,$\QQ(X_{ns}(11))$.

\begin{prop}\label{pull-backs}
The curve $X_{ns}(11)$ admits  morphisms 
defined over $\QQ$ of degree \,$6$, $2$, $2$ and~$6$ to the elliptic curves $A$, $B$, $C$ and $D$ respectively.
Moreover, the corresponding pull-backs of the \mbox{$1$-dimensional} \,$\QQ$-vector spaces of regular differentials
are the $1$-dimensional subspaces of~\,$\Omega^1_{X_{ns}(11)}$ generated by
$$\omega_A={{dX}\over Z}\,,\qquad \omega_B={{dX}\over{2Y+1}}\,,\qquad \omega_C={{dX}\over T}\qquad\hbox{and}\qquad \omega_D={{(3X-1)dX}\over Z}
$$
respectively.
\end{prop}

\noindent{\bf Proof.}\, By Corollary \ref{exotic}, the function field extension
$\QQ(X)\subset \QQ(X,Y,T)$ is Galois, with Galois group isomorphic to Klein's four group.
Since the elliptic curve given by the Weierstrass equation \,$T^2=-(4X^3+7X^2-6X+19)$\, is isomorphic to $C$,
we have the following commutative diagram of degree $2$ morphisms
$$
\xymatrix{
& X_{ns}(11) \ar@{->}[ld]_{\scriptstyle\phi_B} \ar@{->}[d]^{\scriptstyle\phi_H}
\ar@{->}[rd]^{\scriptstyle\phi_C} & \\
B \ar@{->}[rd] & H \ar@{->}[d] & C \ar@{->}[ld] \\
& \PP^1&
}
$$
Here $H$ is the genus $2$ curve given by
$$
Z^2=-(4X^3-4X^2-28X+41)(4X^3+7X^2-6X+19),
$$
and the morphisms \,$\phi_B$, \,$\phi_H$\, and \,$\phi_C$\, are defined as follows:
$$\phi_B(X,Y,T)=(X,Y), \qquad \phi_H(X,Y,T)=\left(X,(2Y+1)T\right)\!, \qquad  \phi_C(X,Y,T)=(X,T).$$
In particular, we can take \,$\omega_B$\, and \,$\omega_C$\, as in the statement.

We now describe degree $6$ morphisms from $X_{ns}(11)$ to the curves $A$ and $D$\, factoring through~$\phi_H$.
To see that $H$ admits degree $3$ morphisms to $A$ and $D$, we use Goursat's formulas as described in the appendix
of~\cite{HOWE}. Substituting \,$X=x+{1\over 3}$\, and \,$Z={{44}\over{3}} z$\, in  the hyperelliptic equation of $H$,
we obtain
$$
t\,z^2\,=\,\left(x^3+3ax+2b\right)\!\left(2dx^3+3cx^2+1\right)
$$
with
$$a=-{{\,22\,}\over 9}\,,\qquad b={{\,847\,}\over{216}}\,,\qquad c={{27}\over{\,242\,}}\,,\qquad d={{9}\over{\,44\,}}
\qquad\hbox{and}\qquad t=-3.\\[3pt]
$$
Note that the discriminants \,$\Delta_1=a^3+b^2$\, and \,$\Delta_2=c^3+d^2$\, are both non-zero.
Then, the maps \,$(x,z)\mapsto(u,v)$,\, with
$$
\begin{array}{ccl}
(u,v) & = & \left(12\,\Delta_1\,\dfrac{-2dx+c}{\,x^3+3ax+2b\,}\,, \ \  z\,\Delta_1\,\dfrac{16dx^3-12cx^2-1}{\,\,(x^3+3ax+2b)^2\,}\right)\!,\\[20pt]
(u,v) & = & \left(12\,\Delta_2\,\dfrac{x^2(ax-2b)}{\,2dx^3+3cx^2+1\,}\,, \ \ z\,\Delta_2\,\dfrac{x^3+12ax-16b}{\,\,(2dx^3+3cx^2+1)^2\,}\right)\!,
\end{array}
$$
are degree $3$ morphisms from $H$ to the genus $1$ curves given by the equations\\[3pt]
$$\begin{array}{ccl}
t\,v^2 & = & u^3 \,+\, 12\left(2a^2d-bc\right)u^2 \,+\, 12\,\Delta_1\left(16ad^2+3c^2\right)u \,+\,
512\,\Delta_1^2\,d^3,\\[10pt]
t\,v^2 & = & u^3 \,+\, 12\left(2bc^2-ad\right)u^2 \,+\, 12\,\Delta_2\left(16b^2c+3a^2\right)u \,+\,
512\,\Delta_2^2\,b^3
\end{array}
$$
respectively. Moreover, the pull-back of the differential \,$du/v$\, of the first curve to $\Omega^1_H$ is
a rational multiple of \,$dx/z$\, and hence of \,$dX/Z$,\,
while the pull-back of the differential \,$du/v$\, of the second curve is a rational multiple
of \,$x\,dx/z$\, and hence of \,$(3X-1)dX/Z$.

Finally, for the above values of \,$a$, $b$, $c$, $d$\, and \,$t$,\, the two genus $1$ curves can be
checked to be isomorphic over $\QQ$ to the elliptic curves $A$ and $D$ respectively. This proves the proposition.

\begin{rem}{\rm  Since the Jacobian of $H$ is isogenous to $A\times D$, we know  that there do exist
non-constant morphisms from $H$ to the curves $A$ and $D$,
but we know of no a priori reason why there should exist morphisms of degree~$3$.
In fact, this was only established by a numerical computation involving the period lattices
of the curves $H$, $A$ and $D$.  Another reason for suspecting that there exist such morphisms
is the fact that the Fourier coefficients of the weight $2$ eigenforms associated to the elliptic
curves $A$ and $D$ are congruent modulo~$3$.}
\end{rem}

\section{\large Automorphisms}\label{aut}

In this section we prove the theorem.  We use the notations of Proposition~\ref{equations} and Proposition~\ref{pull-backs}.

\vspace{1truemm}

Let $\sigma$ be an automorphism of the curve $X_{ns}(11)$. Then $\sigma$ induces an automorphism
of the Jacobian~$J_{ns}(11)$. We recall that this Jacobian is isogenous over~$\QQ$ to the product
of the elliptic curves $A$, $B$, $C$ \,and $D$\, introduced in section~\ref{dif}.

Let us analyze the isogeny relations over $\overline{\QQ}$\, among these four elliptic curves.
The curve~$D$ cannot be isogenous over~$\overline{\QQ}$ to $A$, $B$ or~$C$ because it is the only one
whose $j$-invariant is not integral. The curve $B$ has complex multiplication by the quadratic order of
discriminant~$-11$, so it cannot be isogenous over~$\overline{\QQ}$ to $A$, $C$ or~$D$\, because none of
these three curves admits complex multiplication. Lastly, there is a degree $2$ isogeny between $A$ and $C$
defined over~$\QQ(\sqrt{-11})$.

Therefore, all endomorphisms of \,$J_{ns}(11)$ are defined over $\QQ(\sqrt{-11})$. Furthermore,
the action of \,$\sigma$ on \,$\Omega^1_{X_{ns}(11)}$ with respect to the basis
\,$\omega_B$, $\omega_D$, $\omega_A$, $\omega_C$\, is given by multiplication by a matrix of the form
\begin{equation}
\left(\!\!\begin{array}{rrrr}
\pm 1 & 0 & \ 0 \ & \ 0\\
0 & \pm 1 & \ 0 \ & \ 0\\
0 & 0 & \ a \ & \ b \\
0 & 0 & \ c \ & \ d
\end{array}\!\right)
\end{equation}
for certain \,$a,b,c,d\in\QQ(\sqrt{-11})$. Note that the eigenvalues corresponding to \,$\omega_B$\, and \,$\omega_D$\,
must be roots of unity in this quadratic field, namely $\pm1$, because \,$\sigma$\, has finite order.

\vspace{1truemm}

Let us now consider the functions \,\,$x=\omega_D/\omega_A=3X-1$\, \,and\, \,\,$y=\omega_C/\omega_A=2Y+1$\,
on the elliptic curve $B$. They satisfy the equation\\[-20pt]
\begin{center}
${\frac{1}{4}\,y^2 \,=\, \frac{1}{\,27\,}\,x^3 \,-\, \frac{\,22\,}{9}\,x \,+\, \frac{\,847\,}{108}}\,.$
\end{center}
Then the action of \,$\sigma$ on \,$\Omega^1_{X_{ns}(11)}$ yields
$$
\displaystyle{\sigma(x)= {{\pm x}\over{\,a+cy\,}} \qquad \hbox{and} \qquad \sigma(y)={{\,b+dy\,}\over{\,a+cy\,}}}\,.
$$
In other words, $\sigma$ induces an automorphism of the curve $B$ which, in projective coordinates, is given by
$$
(x\,:\,y\,:\,z)\,\,\longmapsto\,\,(\pm x\,:\,bz + dy\,:\,az + c y).
$$
In particular, $\sigma$ maps the origin $(0:1:0)$ of the elliptic curve $B$ to the point $(0:d:c)$.
This implies $c=0$. Otherwise, the above equation would entail the relation \,$(d/c)^2=847/27$\,
with \,$d/c\in\QQ(\sqrt{-11})$,\, which is impossible. Since the only automorphisms of $B$ fixing
the origin are the identity and the elliptic involution, it follows \,$\sigma(x)=x$\, and \,$\sigma(y)=\pm y$.
Thus, \,$\sigma(X)=X$\, whereas \,$\sigma(Y)$ must be either $Y$ or $1-Y$\!.\,
The equations given for $X_{ns}(11)$ in Proposition~\ref{equations} imply then \,$\sigma(T)=\pm T$\!.\,
This proves the theorem.

\newpage




\vfill


\begin{footnotesize}
\begin{tabular}{ccc}
\begin{tabular}{l}
Valerio Dose \\
\texttt{dose\,\footnotesize{$@$}\,mat.uniroma2.it}\\[5pt]
Ren\'e Schoof\\
\texttt{schoof\,\footnotesize{$@$}\,mat.uniroma2.it}\\[5pt]
Dipartimento di Matematica\\
Universit\`a degli Studi di Roma Tor Vergata\\
Via della Ricerca Scientifica 1\\
00133 Roma, Italy
\end{tabular}

& \qquad\quad &

\begin{tabular}{l}
Josep Gonz\'alez\\
\texttt{josepg\,\footnotesize{$@$}\,ma4.upc.edu}\\[5pt]
Julio Fern\'andez\\
\texttt{julio\,\footnotesize{$@$}\,ma4.upc.edu}\\[5pt]
Departament de Matem\`atica Aplicada 4  \\
Universitat Polit\`ecnica de Catalunya  \\
EPSEVG, Avinguda V\'ictor Balaguer 1\\
08800 Vilanova i la Geltr\'u, Spain
\end{tabular}
\end{tabular}
\end{footnotesize}

\end{document}